\documentclass[a4paper,12pt]{article}
\usepackage{amsmath}
\usepackage{latexsym}
\usepackage{url}
\usepackage{color}
\numberwithin{equation}{section}

\def\R{{\bf R}}

\def\N{{\bf N}}

\def\d{\displaystyle}
\def\e{{\varepsilon}}

\def\v#1{\mbox{\boldmath $#1$}}

\newtheorem{thm}{Theorem}[section]

\newtheorem{prop}{Proposition}[section]
\newtheorem{rem}{Remark}[section]

\title{The lifespan of classical solutions of semilinear
wave equations with spatial weights
and compactly supported data\\
in one space dimension\\
 {\small\it In memory of Professor Masaki Kurokiba}}
\author{
Shunsuke Kitamura
\footnote{
Master course, Mathematical Institute,
Tohoku University,
Aoba, Sendai 980-8578, Japan.
email: shunsuke.kitamura.s8@dc.tohoku.ac.jp (Kitamura),\newline
katsuaki.morisawa.q8@dc.tohoku.ac.jp (Morisawa)},
Katsuaki Morisawa
\footnotemark[1],
Hiroyuki Takamura
\footnote{Mathematical Institute/
Research Alliance Center of Mathematical Sciences,
Tohoku University,
Aoba, Sendai 980-8578, Japan.
e-mail: hiroyuki.takamura.a1@tohoku.ac.jp.}
}
\date{
\[
\begin{array}{ll}
\mbox{\footnotesize{\bf Keywords:}}
& \mbox{\footnotesize semilinear wave equation, one dimension, classical solution, lifespan}\\
\mbox{\footnotesize{\bf MSC2020:}}
& \mbox{\footnotesize primary 35L71, secondary 35B44}\\
\end{array}
\]
}
\begin{document}
\maketitle
\begin{abstract}
This paper studies initial value problems for semilinear wave equations
with spatial weights in one space dimension.
The lifespan estimates of classical solutions for compactly supported data
are established in all the cases of polynomial weights.
The results are classified into two cases according to the total integral of the initial speed.
\end{abstract}

%%%%%%%%%%%%%%%%%%%%%%%%%%%%%%%%%%%%%%%%%%%%%%%%%%
%%%%%%%%%%%%%%%%%%%%% SECTION1 %%%%%%%%%%%%%%%%%%%%%%%
%%%%%%%%%%%%%%%%%%%%%%%%%%%%%%%%%%%%%%%%%%%%%%%%%%

\section{Introduction}
\par
We consider the following initial value problem
for semilinear wave equations with spatial weights. 
\begin{equation}
\label{IVP}
\left\{
\begin{array}{ll}
	\d u_{tt}-\Delta u=\frac{|u|^p}{(1+x^2)^{(1+a)/2}}
	&\mbox{in}\quad \R\times[0,\infty),\\
	u(x,0)=\e f(x),\ u_t(x,0)=\e g(x),
	& x\in\R,
\end{array}
\right.
\end{equation}
where $p>1$, $a\in\R$, $f$ and $g$ are given smooth functions of compact support
and a parameter $\e>0$ is \lq\lq small enough".
\par
When $a=-1$, (\ref{IVP}) is well-studied as a model to ensure 
the optimality of the general theory for nonlinear wave equations.
See Introduction in Imai, Kato, Takamura and Wakasa \cite{IKTW19}
for all the references to this direction including higher dimensions.
More precisely, since we have no time decay of the solution of the free wave equation in one space dimension,
there is no possibility to construct a global-in-time solution of (\ref{IVP}) for any $p>1$.
In fact, we have the finite-time blow-up result by Kato \cite{Kato80}.
Therefore we are interested in the so-called lifespan estimates,
namely, some kind of a stability of a zero solution
because we have an uniqueness of the solution of (\ref{IVP}).
Let $T(\e)$ be, the so-called lifespan, the maximal existence time of the classical solution of (\ref{IVP})
with arbitrary fixed non-zero data.
Due to Zhou \cite{Zhou92}, we have
\begin{equation}
\label{lifespan_zhou}
T(\e)\sim
\left\{
\begin{array}{ll}
C\e^{-(p-1)/2} & \mbox{if}\ \d\int_\R g(x)dx\neq0,\\
C\e^{-p(p-1)/(p+1)} & \mbox{if}\ \d\int_\R g(x)dx=0,
\end{array}
\right.
\end{equation}
where $T(\e)\sim A(\e,C)$ stands for the fact that there are positive constants,
$C_1$ and $C_2$, independent of $\e$ satisfying $A(\e,C_1)\le T(\e)\le A(\e,C_2)$.
We note that $p>1$ implies
\[
\frac{p-1}{2}<\frac{p(p-1)}{p+1},
\]
so that the first quantity is smaller than the second one in (\ref{lifespan_zhou}).
This phenomenon follows from the fact that Huygens' principle holds
if the total integral of the initial speed is zero.

\par
When $a\neq-1$, there are a few results only with the assumption that
the data has non-compact support.
This kind of the problem was first proposed by Suzuki \cite{Suzuki10}
in which the nonlinearity $|u|^p$ is replaced by $|u|^{p-1}u$ 
showing the global-in-time existence for odd function data
when $p>(1+\sqrt{5})/2$ and $pa>1$.
She also studied the blow-up result of modified integral equations.
See the section 6 in \cite{Suzuki10}.
Later, Kubo, Osaka and Yazici \cite{KOY13} extended such a result for all $p>1$ and $pa>1$.
Moreover, they obtained  the blow-up in finite-time for (\ref{IVP}) with some positive data
for $p>1$ and $a\ge-1$.
Inspired by some computation of the upper bound of the lifespan in \cite{KOY13},
Wakasa \cite{Wakasa17} obtained the following lifespan estimate for (\ref{IVP}).
\begin{equation}
\label{lifespan_wakasa}
T(\e)\sim
\left\{
\begin{array}{ll}
C\e^{-(p-1)/(1-a)} & \mbox{for}\ -1\le a<0,\\
\phi^{-1}(C\e^{-(p-1)}) & \mbox{for}\ a=0,\\
C\e^{-(p-1)} & \mbox{for}\ a>0,
\end{array}
\right.
\end{equation}
where $\phi^{-1}$ is an inverse function of $\phi$ defined by
\begin{equation}
\label{phi}
\phi(s):=s\log(2+s).
\end{equation}
We note that this result is also available even if $|u|^p$ is replaced with $|u|^{p-1}u$ in (\ref{IVP}).

\par
The aim of this paper is to establish the lifespan estimates
for compactly supported data in all the cases of $a$ including $a<-1$. 
More precisely, our results are the following estimates.
\begin{equation}
\label{lifespan_non-zero}
T(\e)\sim
\left\{
\begin{array}{ll}
C\e^{-(p-1)/(1-a)} & \mbox{for}\ a<0,\\
\phi^{-1}(C\e^{-(p-1)}) & \mbox{for}\ a=0,\\
C\e^{-(p-1)} & \mbox{for}\ a>0
\end{array}
\right.
\quad
\mbox{if}\ \int_\R g(x)dx\neq0
\end{equation}
and
\begin{equation}
\label{lifespan_zero}
T(\e)\sim
\left\{
\begin{array}{ll}
C\e^{-p(p-1)/(1-pa)} & \mbox{for}\ a<0,\\
\psi_p^{-1}(C\e^{-p(p-1)}) & \mbox{for}\ a=0,\\
C\e^{-p(p-1)} & \mbox{for}\ a>0
\end{array}
\right.
\quad
\mbox{if}\ \int_\R g(x)dx=0,
\end{equation}
where $\psi_p^{-1}$ is an inverse function of $\psi_p$ defined by 
\begin{equation}
\label{psi}
\psi_p(s):=s\log^p(2+s).
\end{equation}
We remark that the quantities in all the cases of (\ref{lifespan_non-zero}) are larger than
those of (\ref{lifespan_zero}).
This fact follows from the trivial inequality
\[
\frac{p-1}{1-a}<\frac{p(p-1)}{1-pa}
\]
by $p>1$ in the first case of $a<0$.
For the second case $a=0$, one can check it by comparing two functions
$\phi^{-1}(\xi)$ and $\psi_p^{-1}(\xi^p)$ with respect to the large variable $\xi$
by making use of differentiation.
The third case $a>0$ is trivial.
We also note that Suzuki obtained
$T(\e)<\infty$ for $-1\le a\le 1$ and $g(x)\ge0(\not \equiv0)$ in the section 7 in \cite{Suzuki10}.
Her original result is established for the nonlinear term $|u|^{p-1}u$,
but the proof of $u\ge0$ is missing for compactly supported data.

\par
It is interesting to compare the nonlinear term with time-decaying weights
in Kato, Takamura and Wakasa \cite{KTW19}
which is closely related to the scale-invariantly damped wave equations.
In such a situation, we have a possibility to obtain the global-in-time existence
for the super-critical case,
and the exponential type estimate of the lifespan for the critical case.

\par
This paper is organized as follows.
In the next section, (\ref{lifespan_non-zero}) and (\ref{lifespan_zero}) are divided into four theorems,
and the preliminaries for their proofs are introduced.
Section 3 and 4 are devoted to the proofs of the longtime existence
and the blow-up in finite time of the solution, respectively.
The main method in this paper is based on point-wise estimates 
which are originally introduced by John \cite{John79} in three space dimensions,
and developed by Zhou \cite{Zhou92} and Kato, Takamura and Wakasa \cite{KTW19}
in one space dimension.

%%%%%%%%%%%%%%%%%%%%%%%%%%%%%%%%%%%%%%%%%%%%%%%%%%
%%%%%%%%%%%%%%%%%%%%% SECTION2 %%%%%%%%%%%%%%%%%%%%%%%
%%%%%%%%%%%%%%%%%%%%%%%%%%%%%%%%%%%%%%%%%%%%%%%%%%

\section{Main results and preliminaries}

Throughout of this paper, we assume that the initial data
$(f,g)\in C_0^2(\R)\times C^1_0(\R)$ satisfies
\begin{equation}
\label{supp_initial}
\mbox{\rm supp}\  (f,g)\subset\{x\in\R:|x|\le R\},\quad R\ge1.
\end{equation}
Our results on (\ref{lifespan_non-zero}) and (\ref{lifespan_zero})
are splitted into the following four theorems.

\begin{thm}
\label{thm:lower-bound_non-zero}
Assume the support condition (\ref{supp_initial})
and
\begin{equation}
\label{initial_non-zero}
\int_\R g(x)dx\neq0.
\end{equation}
Then, there exists a positive constant $\e_1=\e_1(f,g,p,a,R)>0$ such that
a classical solution $u\in C^2(\R\times[0,T))$ of (\ref{IVP}) exists as far as $T$ satisfies
\begin{equation}
\label{lower-bound_non-zero}
T\le
\left\{
\begin{array}{ll}
c\e^{-(p-1)/(1-a)} & \mbox{for}\ a<0,\\
\phi^{-1}(c\e^{-(p-1)}) & \mbox{for}\ a=0,\\
c\e^{-(p-1)} & \mbox{for}\ a>0,
\end{array}
\right.
\end{equation}
where $0<\e\le\e_1$, $c$ is a positive constant independent of $\e$
and $\phi$ is the one in (\ref{phi}).
\end{thm}

\begin{rem}
\label{rem:wakasa}
In Wakasa \cite{Wakasa17} for the non-compactly supported data,
the assumption on the data is
\[
f\in C^2(\R)\mbox{ with }\|f\|_{L^\infty(\R)}<\infty,
\quad
g\in C^1(\R)\mbox{ with }\|g\|_{L^1(\R)}<\infty
\]
without (\ref{initial_non-zero}),
so that the case of $a\ge-1$ in Theorem \ref{thm:lower-bound_non-zero}
is already established by (\ref{lifespan_wakasa}).
\end{rem}

\begin{thm}
\label{thm:lower-bound_zero}
Assume the support condition (\ref{supp_initial})
and
\begin{equation}
\label{initial_zero}
\int_\R g(x)dx=0.
\end{equation}
Then, there exists a positive constant $\e_2=\e_2(f,g,p,a,R)>0$ such that
a classical solution $u\in C^2(\R\times[0,T))$ of (\ref{IVP}) exists as far as $T$ satisfies
\begin{equation}
\label{lower-bound_zero}
T\le
\left\{
\begin{array}{ll}
c\e^{-p(p-1)/(1-pa)} & \mbox{for}\ a<0,\\
\psi_p^{-1}(c\e^{-p(p-1)}) & \mbox{for}\ a=0,\\
c\e^{-p(p-1)} & \mbox{for}\ a>0,
\end{array}
\right.
\end{equation}
where $0<\e\le\e_2$, $c$ is a positive constant independent of $\e$ and
$\psi_p$ is the one in (\ref{psi}).
\end{thm}

\begin{thm}
\label{thm:upper-bound_non-zero}
Assume the support condition (\ref{supp_initial})
and
\begin{equation}
\label{positive_non-zero}
\int_{\R}g(x)>0.
\end{equation}
Then, there exists a positive constant $\e_3=\e_3(g,p,a,R)>0$ such that
a classical solution $u\in C^2(\R\times[0,T))$ of (\ref{IVP}) cannot exist whenever $T$ satisfies
\begin{equation}
\label{upper-bound_non-zero}
T\ge
\left\{
\begin{array}{ll}
C\e^{-(p-1)/(1-a)} & \mbox{for}\ a<0,\\
\phi^{-1}(C\e^{-(p-1)}) & \mbox{for}\ a=0,\\
C\e^{-(p-1)} & \mbox{for}\ a>0,
\end{array}
\right.
\end{equation}
where $0<\e\le\e_3$, $C$ is a positive constant independent of $\e$
and $\phi$ is the one in (\ref{phi}).
\end{thm}

\begin{thm}
\label{thm:upper-bound_zero}
Assume the support condition (\ref{supp_initial})
and
\begin{equation}
\label{positive_zero}
f(x)\ge0(\not\equiv0),\quad g(x)\equiv0.
\end{equation}
Then, there exists a positive constant $\e_4=\e_4(f,p,a,R)>0$ such that
a classical solution $u\in C^2(\R\times[0,T))$ of (\ref{IVP}) cannot exist whenever $T$ satisfies
\begin{equation}
\label{upper-bound_zero}
T\ge
\left\{
\begin{array}{ll}
C\e^{-p(p-1)/(1-pa)} & \mbox{for}\ a<0,\\
\psi_p^{-1}(C\e^{-p(p-1)}) & \mbox{for}\ a=0,\\
C\e^{-p(p-1)} & \mbox{for}\ a>0,
\end{array}
\right.
\end{equation}
where $0<\e\le\e_4$, $C$ is a positive constant independent of $\e$
and $\psi_p$ is the one in (\ref{psi}).
\end{thm}

\begin{rem}
In view of the definition of lifespan $T(\e)$, 
Theorems \ref{thm:lower-bound_non-zero} and \ref{thm:upper-bound_non-zero}
imply (\ref{lifespan_non-zero}), also
Theorems  \ref{thm:lower-bound_zero} and \ref{thm:upper-bound_zero}
imply (\ref{lifespan_zero}).
\end{rem}

\par
All the proofs of above theorems are given in following sections.
Here we shall introduce preliminaries.
Let $u$ be a classical solution of (\ref{IVP}) in the time interval $[0,T)$.
Then the support condition of the initial data, (\ref{supp_initial}), implies that
\begin{equation}
\label{support_sol}
\mbox{supp}\ u(x,t)\subset\{(x,t)\in\R\times[0,T):|x|\le t+R\}.
\end{equation}
For example, see Appendix of John \cite{John_book} for this fact.
It is well-known that $u$ satisfies the following integral equation.
\begin{equation}
\label{integral}
u(x,t)=\e u^0(x,t)+L_a(|u|^p)(x,t),
\end{equation}
where $u^0$ is a solution of the free wave equation with the same initial data,
\begin{equation}
\label{linear}
u^0(x,t):=\frac{1}{2}\{f(x+t)+f(x-t)\}+\frac{1}{2}\int_{x-t}^{x+t}g(y)dy,
\end{equation}
and a linear integral operator $L_a$ for a function $v=v(x,t)$ in Duhamel's term is defined by
\begin{equation}
\label{nonlinear}
L_a(v)(x,t):=\frac{1}{2}\int_0^tds\int_{x-t+s}^{x+t-s}\frac{v(y,s)}{(1+y^2)^{(1+a)/2}}dy.
\end{equation}

\begin{prop}
\label{prop:continuity}
Assume that $(f,g)\in C^2(\R)\times C^1(\R)$.
Let $u$ be a continuous solution of (\ref{integral}). 
Then, $u$ is a classical solution of (\ref{IVP}).
\end{prop}

\par\noindent
{\bf Proof.} In view of (\ref{nonlinear}), 
the differentiability of $L_a(v)$ follows from the continuity of $v$.
Therefore the conclusion follows from the regularity assumption on the initial data.
\hfill$\Box$

\vskip10pt
\par
The following property, namely Huygens' principle, of $u^0$ will play an essential role
in the proofs of Theorems \ref{thm:lower-bound_zero} and \ref{thm:upper-bound_zero}. 

\begin{prop}
\label{prop:huygens}
Assume (\ref{supp_initial}) and (\ref{initial_zero}).
Then, $u^0$ in (\ref{linear}) satisfies
\begin{equation}
\label{huygens}
\mbox{\rm supp}\ u^0(x,t)\subset\{(x,t)\in\R\times[0,\infty):(t-R)_+\le|x|\le t+R\}.
\end{equation}
\end{prop}

\par\noindent
{\bf Proof.} For $t\ge R$ and $|x|\le t-R$, we have
\[
x+t\ge R\quad\mbox{and}\quad x-t\le-R.
\]
Therefore it follows from (\ref{supp_initial}), (\ref{initial_zero}) and (\ref{linear}) that
\[
u^0(x,t)\equiv0\quad\mbox{for}\ t\ge R\quad\mbox{and}\quad |x|\le t-R.
\]
On the other hand, it is trivial that
\[
u^0(x,t)\equiv0\quad\mbox{for}\ t+R<|x|,
\]
so that (\ref{huygens}) holds.
\hfill$\Box$

\vskip10pt
\par
Due to Proposition \ref{prop:huygens} as well as (\ref{support_sol}),
we shall divide the support of the solution into three pieces,
the interior domain
\begin{equation}
\label{Dint}
D_{\rm Int}:=\{(x,t)\in\R\times[0,T]:t+|x|\ge R,\ t-|x|\ge R\},
\end{equation}
the exterior domain
\begin{equation}
\label{Dext}
D_{\rm Ext}:=\{(x,t)\in\R\times[0,T]:t+|x|\ge R,\ \left|t-|x|\right|\le R\},
\end{equation}
and the small domain near the origin
\begin{equation}
\label{Dori}
D_{\rm Ori}:=\{(x,t)\in\R\times[0,T]:t+|x|\le R,\ \left|t-|x|\right|\le R\}.
\end{equation}
We will see that the lifespan is determined by point-wise estimates of the solution in $D_{\rm Int}$.

%%%%%%%%%%%%%%%%%%%%%%%%%%%%%%%%%%%%%%%%%%%%%%%%%%
%%%%%%%%%%%%%%%%%%%%% SECTION3 %%%%%%%%%%%%%%%%%%%%%%%
%%%%%%%%%%%%%%%%%%%%%%%%%%%%%%%%%%%%%%%%%%%%%%%%%%

\section{Proofs of Theorems \ref{thm:lower-bound_non-zero} and \ref{thm:lower-bound_zero}}

In this section, we investigate the lower bound of the lifespan.
In view of Remark \ref{rem:wakasa},
only the case of $a<-1$ should be considered in the proof of Theorem \ref{thm:lower-bound_non-zero}.
But, following the proof of Wakasa \cite{Wakasa17},
all the estimates for the case of $-1\le a<0$ hold also for the case of $a<-1$,
so that we can omit its proof here.
In fact, we have to show that (4.6) in Wakasa \cite{Wakasa17}
is also established for $a<-1$.
The case of $0\le x\le t\le T$ is trivial and another case of $x\ge t$ follows from
$|x|\le t+R$ by (\ref{support_sol}). 
\par
From now on, we shall prove Theorem \ref{thm:lower-bound_zero} only.
To this end, we have to set the following function space which is different from Wakasa \cite{Wakasa17}.
Following Kato, Takamura and Wakasa \cite{KTW19},
we shall construct a solution as a limit of the sequence $\{U_n(x,t)\}_{n\in\N}$ defined by
\begin{equation}
\label{sequence}
U_{n+1}=L_a(|U_n+\e u^0|^p),\quad U_1\equiv0
\end{equation}
in the weighted $L^\infty$ space.
Let $w$ a weight function defined by
\begin{equation}
w(r,t):=
\left\{
\begin{array}{ll}
(t+r+3R)^a & \mbox{for}\ a<0,\\
\{\log(t+r+3R)\}^{-1} & \mbox{for}\ a=0,\\
1 & \mbox{for}\ a>0
\end{array}
\right.
\end{equation}
and a weighted norm of a function $U=U(x,t)$ by
\begin{equation}
\label{norm}
\|U\|:=\sup_{(x,t)\in\R\times[0,T]}w(|x|,t)|U(x,t)|.
\end{equation}
We note that H\"older's inequality
\begin{equation}
\label{Holder}
\||U|^\theta|V|^{1-\theta}\|\le\|U\|^\theta\|V\|^{1-\theta}\qquad(0\le\theta\le1)
\end{equation}
holds.

\par
Then we have a priori estimates in the following propositions.

\begin{prop}
\label{prop:apriori_0}
Suppose that the assumption of Theorem \ref{thm:lower-bound_zero} is fulfilled.
Let $L_a$ and $u^0$ be the ones in (\ref{linear}) and (\ref{nonlinear}) respectively. 
Then, for $U=U(x,t)\in C^0_0(\R\times[0,T])$ with supp\ $U\subset\{(x,t)\in\R\times[0,T]:|x|\le t+R\}$,
there exists a positive constant $M=M(f,g,a,m,R)$ such that
\begin{equation}
\label{linear-bound}
\|L_a(|u^0|^{p-m}|U|^m)\|\le M\{\|U\|D(T)\}^m
\quad\mbox{for}\ m=0,1,
\end{equation}
where $D(T)$ is defined by
\begin{equation}
\label{D(T)}
D(T):=
\left\{
\begin{array}{ll}
(T+2R)^{-a} & \mbox{for}\ a<0,\\
\log(T+3R) & \mbox{for}\  a=0,\\
1 & \mbox{for}\ a>0.
\end{array}
\right.
\end{equation}
\end{prop}

\begin{prop}
\label{prop:apriori}
Suppose that the assumption of Theorem \ref{thm:lower-bound_zero} is fulfilled.
Let $L_a$ be the one in (\ref{nonlinear}). 
Then, for $U=U(x,t)\in C^0_0(\R\times[0,T])$ with supp\ $U\subset\{(x,t)\in\R\times[0,T]:|x|\le t+R\}$,
 there exists a positive constant $C=C(f,g,a,R)$ such that
\begin{equation}
\|L_a(|U|^p)\|\le C\|U\|^pE(T),
\end{equation}
where $E(T)$ is defiend by
\begin{equation}
\label{E}
E(T):=
\left\{
\begin{array}{ll}
(T+2R)^{1-pa} & \mbox{for}\ a<0,\\
(T+R)\log^p(T+3R) & \mbox{for}\  a=0,\\
T+R & \mbox{for}\ a>0.
\end{array}
\right.
\end{equation}
\end{prop}

%%%%%%%%%%%%%%%%%%%%%%%%%%
\par
First we shall prove the main theorem.
The proofs of the propositions above are given later.

\vskip10pt
\par\noindent
{\bf Proof of Theorem \ref{thm:lower-bound_zero}.}
\par
By virtue of Proposition \ref{prop:continuity},
it is sufficient to construct a continuous solution of the integral equation (\ref{integral}). 
Following Kato, Takamura and Wakasa \cite{KTW19}, let $X$ be a Banach space defined by
\[
X:=\{ U(x,t) \in C(\R \times [0,T]) : \mbox{supp} \ U\subset \{(x,t)\in\R\times[0,T]:|x| \leq t+R \} \}
\]
which is equipped with the norm (\ref{norm}), and its closed subspace $Y$ by
\[
Y:=\{U\in X:\|U\|\le2^{p+1}M\e^p\},
\]
where $M$ is the one in Proposition \ref{prop:apriori_0}.
We note that $\{U_n\}$ in (\ref{sequence}) is the sequence in $X$ because
\[
\mbox{supp}\ u^0\subset\{(x,t):|x|\le t+R\}
\]
is trivial and
\[
\mbox{supp}\ U_n\subset\{(y,s):|y|\le s+R \}
\]
implies that $(x,t)\in\mbox{supp}\ U_{n+1}$ satisfies
\[
|x|\le |y|+t-s\le t+R
\]
because
\[
x-t+s\le y\le x+t-s
\]
for $y$ in the domain of the integral in $L_a$ is equivalent to
\[
|y-x|\le t-s.
\]
The continuity of the sequence is also trivial.

\par
Since we have
\[
|U_{n+1}|\le2^p\{L_a(|U_n|^p)+\e^pL_a(|u^0|^p)\},
\]
Propositions \ref{prop:apriori_0} with $m=0$ and \ref{prop:apriori} yield
\[
\|U_{n+1}\|\le2^pC\|U_n\|^pE(T)+2^pM\e^p,
\]
where $C$ is the one in Proposition \ref{prop:apriori}.
Hence the boundedness in $Y$ of $\{U_n\}$ in (\ref{sequence}) follows from
\begin{equation}
\label{convergence1}
2^{p^2+2p}CM^pE(T)\e^{p^2}\le2^pM\e^p.
\end{equation}

\par
From now on, we assume (\ref{convergence1}).
Since
\[
\begin{array}{ll}
|U_{n+1}-U_n|
&\le L_a\left(\left||U_n+\e u^0|^p-|U_{n-1}+\e u^0|^p\right|\right)\\
&\le pL_a\left(|U_{n-1}+\e u^0+\theta(U_n-U_{n-1})|^{p-1}|U_n-U_{n-1}|\right)\\
&\le3^{p-1}pL_a\{(|U_n|^{p-1}+|U_{n-1}|^{p-1}+\e^{p-1}|u^0|^{p-1})|U_n-U_{n-1}|\}
\end{array}
\]
holds with some $\theta\in(0,1)$,
Propositions \ref{prop:apriori_0} with $m=1$ and \ref{prop:apriori} yield
\[
\begin{array}{ll}
\|U_{n+1}-U_n\|\le
&3^{p-1}pC(\|U_n\|^{p-1}+\|U_{n-1}\|^{p-1})\|U_n-U_{n-1}\|E(T)\\
&+3^{p-1}p\e^{p-1}M\|U_n-U_{n-1}\|D(T).
\end{array}
\]
Here we have employed (\ref{Holder}) as
\[
\begin{array}{ll}
\|L_a(|U_n|^{p-1}|U_n-U_{n-1}|)\|
&=\|L_a\{(|U_n|^{1-1/p}|U_n-U_{n-1}|^{1/p})^p\}\|\\
&\le C\||U_n|^{1-1/p}|U_n-U_{n-1}|^{1/p}\|^pE(T)\\
&\le C\|U_n\|^{p-1}\|U_n-U_{n-1}\|E(T)
\end{array}
\]
and so on.

\par
Hence $\{U_m\}$ is a Cauchy sequence in $Y$ provided
\begin{equation}
\label{convergence2}
3^{p-1}pC\cdot2(2^{p+1}M\e^p)^{p-1}E(T)+3^{p-1}pM\e^{p-1}D(T)\le\frac{1}{2}.
\end{equation} 
We note that
(\ref{convergence1}) and (\ref{convergence2}) guarantee the existence of a limit of $\{U_n\}$ in $Y$.

\par
When $a>0$, it is easy to find $c$ and $\e_2$ in (\ref{lower-bound_zero})
because of $D(T)=1$ and $E(T)=T+R$.
We omit details.
\par
When $a=0$, let us look for a sufficient condition on $T$ to (\ref{convergence1}) and (\ref{convergence2}).
The definitions of $D(T)$ and $E(T)$ in (\ref{D(T)}) and (\ref{E}) respectively yield
\begin{equation}
\label{convergence_a=0}
\left\{
\begin{array}{l}
2^{p^2+p}CM^{p-1}\e^{p(p-1)}(T+R)\log^p(T+3R)\le1,\\
2^{p^2+1}3^{p-1}pCM^{p-1}\e^{p(p-1)}(T+R)\log^p(T+3R)\\
\qquad+2\cdot3^{p-1}pM\e^{p-1}\log(T+3R)\le1.
\end{array}
\right.
\end{equation}
Assume that
\[
T\ge R.
\]
Then (\ref{convergence_a=0}) follows from
\[
\left\{
\begin{array}{l}
2^{p^2+2p+1}CM^{p-1}\e^{p(p-1)}T\log^p(T+2)\le1,\\
2^{p^2+p+3}3^{p-1}pCM^{p-1}\e^{p(p-1)}T\log^p(T+2)\le1,\\
2^33^{p-1}pM\e^{p-1}\log(T+2)\le1
\end{array}
\right.
\]
because of 
\[
\log(T+3R)\le2\log(T+2)\quad \mbox{for}\ T\ge R.
\]
Therefore Theorem \ref{thm:lower-bound_zero} for $a=0$ is established with
\[
T\le\psi_p^{-1}(C'\e^{-p(p-1)})
\quad\mbox{for}\ 0<\e\le\e_2,
\]
where
\[
C':=\left(2^{p^2+2p+2}3^{p-1}pCM^{p-1}\right)^{-1}>0
\]
and a number $\e_2$ is defined to satisfy
\[
R\le\psi_p^{-1}(C'\e_2^{-p(p-1)})\le\exp(2^{-3}3^{1-p}p^{-1}M^{-1}\e_2^{-(p-1)})-2.
\]
This is possible. The first inequality is trivial. Setting
\[
\Psi(s):=\exp(2^{-3}3^{1-p}p^{-1}M^{-1}s)-2-\psi_p^{-1}(C's^p),
\]
we have
\[
\Psi'(s)=2^{-3}3^{1-p}p^{-1}M^{-1}\exp(2^{-3}3^{1-p}p^{-1}M^{-1}s)-\frac{pC's^{p-1}}{\psi_p'(C's^p)}
\]
where
\[
\psi_p'(s)=\log^p(2+s)+\frac{ps}{2+s}\log^{p-1}(2+s).
\]
Hence the second inequality can be valid by taking $\e_3^{-(p-1)}$ large enough
because it is easy to find a point $s_0$ independent of $\e$ such that
\[
\Psi'(s)\ge1\quad\mbox{for}\ s\ge s_0.
\]

\par
The case of $a<0$ is almost similar to the above.
(\ref{convergence1}) and (\ref{convergence2}) follow from
\[
\left\{
\begin{array}{l}
2^{p^2+p}CM^{p-1}\e^{p(p-1)}(T+2R)^{1-pa}\le1,\\
2^{p^2+1}3^{p-1}pCM^{p-1}\e^{p(p-1)}(T+2R)^{1-pa}\\
\qquad+2\cdot3^{p-1}pM\e^{p-1}(T+2R)^{-a}\le1.
\end{array}
\right.
\]
Since
\[
\frac{p(p-1)}{1-pa}\le\frac{p-1}{-a}
\]
holds, it is easy to see that (\ref{lower-bound_zero}) for $a<0$ is established.
Therefore the proof of Theorem \ref{thm:lower-bound_zero} is now completed.
\hfill$\Box$

%%%%%%%%%%%%%%%%%%%%%%%%%%
\vskip10pt
\par\noindent
{\bf Proof of Proposition \ref{prop:apriori_0}.}
\par
In view of Proposition \ref{prop:huygens} and (\ref{linear}), we have
\[
|L_a(|u^0|^{p-m}|U|^m)(x,t)|\le\frac{(C_{f,g})^{p-m}\|U\|^m}{2}I_0(x,t),
\]
where
\begin{equation}
\label{I_0}
I_0(x,t):=\int_0^tds
\int_{x-t+s}^{x+t-s}\frac{w(|y|,s)^{-m}\chi_0(y,s)}{(1+y^2)^{(1+a)/2}}dy,
\end{equation}
\begin{equation}
\label{chi_0}
\chi_0(y,s):=
\left\{
\begin{array}{ll}
1 & \mbox{for}\ (s-R)_+\le|y|\le s+R\\
0 & \mbox{otherwise}
\end{array}
\right.
\end{equation}
and
\[
C_{f,g}:=\|f\|_{L^\infty(\R)}+\frac{1}{2}\|g\|_{L^1(\R)}>0.
\]
Therefore Proposition \ref{prop:apriori_0} follows from
\begin{equation}
\label{basic-estimate_0}
I_0(x,t)\le Mw(|x|,t)^{-1}D(T)^m
\quad\mbox{for}\ (x,t)\in D_{\rm Ext}\cup D_{\rm Ori}.
\end{equation}
Due to the symmetry of $I_0$ on $x$ as $I_0(-x,t)=I_0(x,t)$,
it is sufficient to show (\ref{basic-estimate_0}) in case of 
\[
x\ge0.
\]

\par
From now on, all the constants $C=C(f,g,a,m,R)>0$ may change from line to line for simplicity.
Changing variables by
\begin{equation}
\label{coordinate}
\alpha:=s+y,\quad\beta:=s-y
\end{equation}
and making use of
\begin{equation}
\label{equivalent}
\frac{1}{\sqrt{2}}(1+|y|)\le\sqrt{1+y^2}\le1+|y|
\quad\mbox{for}\ y\in\R,
\end{equation}
we have that
\[
I_0(x,t)\le C
\left\{
\begin{array}{ll}
I_{01}(x,t)+I_{02}(x,t) & \mbox{for}\ (x,t)\in D_{\rm Ext},\\
I_{03}(x,t) &\mbox{for}\ (x,t)\in D_{\rm Ori},
\end{array}
\right.
\]
where
\[
\begin{array}{l}
\d I_{01}(x,t):=\int_{-R}^{t-x}d\beta\int_R^{t+x}\frac{w(y,s)^{-m}}{(1+(\alpha-\beta)/2)^{1+a}}d\alpha,\\
\d I_{02}(x,t):=\int_{-R}^{t-x}d\beta\int_{-\beta}^R\frac{w(|y|,s)^{-m}}{(1+|\alpha-\beta|/2)^{1+a}}d\alpha,\\
\d I_{03}(x,t):=\int_{-t-x}^{t-x}d\beta\int_{-\beta}^{t+x}\frac{w(|y|,s)^{-m}}{(1+|\alpha-\beta|/2)^{1+a}}d\alpha.
\end{array}
\]

\par
First, we shall estimate $I_{01}$ in $D_{\rm Ext}$.
Extending the domain of the integral, we have
\[
I_{01}(x,t)\le \int_{-R}^Rd\beta\int_R^{t+x}\frac{w(y,s)^{-m}}{(1+(\alpha-\beta)/2)^{1+a}}d\alpha.
\]
When $a>0$, the $\alpha$-integral is dominated by
\[
\left[\frac{2}{-a}\left(1+\frac{\alpha-\beta}{2}\right)^{-a}\right]_{\alpha=R}^{\alpha=t+x}
\le\frac{2}{a}.
\]
When $a=0$, the $\alpha$-integral is dominated by
\[
\begin{array}{l}
\d\log^m(t+x+3R)\left[2\log\left(1+\frac{\alpha-\beta}{2}\right)\right]_{\alpha=R}^{\alpha=t+x}\\
\d\le2^{1+m}\log(t+x+3R)\log^m(T+3R)
\end{array}
\]
because of
\[
\log(t+x+3R)\le\log(2t+4R)\le2\log(T+3R).
\]
When $a<0$, the $\alpha$-integral is dominated by
\[
\begin{array}{l}
\d(t+x+3R)^{m(-a)}\left[\frac{2}{-a}\left(1+\frac{\alpha-\beta}{2}\right)^{-a}\right]_{\alpha=R}^{\alpha=t+x}\\
\d\le\frac{2^{1+a+m(-a)}}{-a}(t+x+3R)^{-a}(T+2R)^{m(-a)}.
\end{array}
\]
Hence we obtain
\[
I_{01}(x,t)\le Cw(x,t)^{-1}D(T)^m\quad\mbox{for}\ (x,t)\in D_{\rm Ext}.
\]
On the other hand, it is easy to see that
\[
I_{02}(x,t)\le C\int_{-R}^Rd\beta\int_{-R}^R\frac{w(|y|,s)^{-m}}{(1+|\alpha-\beta|/2)^{1+a}}d\alpha\le C
\quad\mbox{for}\ (x,t)\in D_{\rm Ext}.
\]
Moreover, similarly to $I_{02}$ in $D_{\rm Ext}$, we also have
\[
I_{03}(x,t)\le C\int_{-R}^Rd\beta\int_{-R}^R\frac{w(|y|,s)^{-m}}{(1+|\alpha-\beta|/2)^{1+a}}d\alpha\le C
\quad\mbox{for}\ (x,t)\in D_{\rm Ori}.
\]
Therefore, summing up, we obtain (\ref{basic-estimate_0}) as desired.
\hfill$\Box$

%%%%%%%%%%%%%%%%%%%%%%%%%%%%%%%%
\vskip10pt
\par\noindent
{\bf Proof of Proposition \ref{prop:apriori}.}
\par
The proof is almost similar to the one of Proposition \ref{prop:apriori_0}.
Due to (\ref{support_sol}), we have
\[
|L_a(|U|^p)(x,t)|\le\frac{\|U\|^p}{2}\int_0^tds
\int_{x-t+s}^{x+t-s}\frac{w(|y|,s)^{-p}}{(1+y^2)^{(1+a)/2}}\chi(y,s)dy,
\]
where
\begin{equation}
\label{chi}
\chi(y,s):=
\left\{
\begin{array}{ll}
1 & \mbox{for}\ |y|\le s+R,\\
0 & \mbox{otherwise}.
\end{array}
\right.
\end{equation}
Therefore Proposition \ref{prop:apriori} follows from
\begin{equation}
\label{basic-estimate}
I(x,t)\le CE(T)w(|x|,t)^{-1}
\quad\mbox{for}\ (x,t)\in D_{\rm Int}\cup D_{\rm Ext}\cup D_{\rm Ori},
\end{equation}
where
\begin{equation}
\label{I}
I(x,t):=\int_0^tds\int_{x-t+s}^{x+t-s}\frac{w(|y|,s)^{-p}}{(1+y^2)^{(1+a)/2}}\chi(y,s)dy.
\end{equation}
Similarly to $I_0$ in the proof of Proposition \ref{prop:apriori_0},
it is sufficient to show (\ref{basic-estimate}) in case of 
\[
x\ge0.
\]

\par
From now on, all the constants $C=C(f,g,a,R)>0$ may change from line to line for simplicity.
Changing variables by (\ref{coordinate}) again, we have that
\[
I(x,t)\le C
\left\{
\begin{array}{ll}
I_{11}(x,t)+I_{12}(x,t)+I_{13}(x,t)+I_{14}(x,t) & \mbox{for}\ (x,t)\in D_{\rm Int},\\
I_{21}(x,t)+I_{22}(x,t) & \mbox{for}\ (x,t)\in D_{\rm Ext},\\
I_{3}(x,t) &\mbox{for}\ (x,t)\in D_{\rm Ori},
\end{array}
\right.
\]
where
\[
\begin{array}{l}
\d I_{11}(x,t):=\int_R^{t-x}d\beta
\int_R^{t+x}\frac{w(|y|,s)^{-p}}{(1+|\alpha-\beta|/2)^{1+a}}d\alpha,\\
\d I_{12}(x,t):=\int_{-R}^Rd\beta
\int_R^{t+x}\frac{w(y,s)^{-p}}{(1+(\alpha-\beta)/2)^{1+a}}d\alpha,\\
\d I_{13}(x,t):=\int_R^{t-x}d\beta
\int_{-R}^R\frac{w(-y,s)^{-p}}{(1-(\alpha-\beta)/2)^{1+a}}d\alpha,\\
\d I_{14}(x,t):=\int_{-R}^Rd\beta
\int_{-\beta}^R\frac{w(|y|,s)^{-p}}{(1+|\alpha-\beta|/2)^{1+a}}d\alpha,\\
\d I_{21}(x,t):=\int_{-R}^{t-x}d\beta
\int_R^{t+x}\frac{w(y,s)^{-p}}{(1+(\alpha-\beta)/2)^{1+a}}d\alpha,\\
\d I_{22}(x,t):=\int_{-R}^{t-x}d\beta
\int_{-\beta}^R\frac{w(|y|,s)^{-p}}{(1+|\alpha-\beta|/2)^{1+a}}d\alpha,\\
\d I_3(x,t):=\int_{-t-x}^{t-x}d\beta
\int_{-\beta}^{t+x}\frac{w(|y|,s)^{-p}}{(1+|\alpha-\beta|/2)^{1+a}}d\alpha.
\end{array}
\]

\par
First, we shall estimate $I_{11}$ in $D_{\rm Int}$.
Since the symmetry of the integrand in $y=(\alpha-\beta)/2$, we have
\[
\begin{array}{l}
\d\int_R^{t-x}d\beta
\int_R^{t-x}\frac{w(|y|,s)^{-p}}{(1+|\alpha-\beta|/2)^{1+a}}d\alpha
\\
\d=2\int_R^{t-x}d\beta
\int_\beta^{t-x}\frac{w(|y|,s)^{-p}}{(1+|\alpha-\beta|/2)^{1+a}}d\alpha,
\end{array}
\]
so that we obtain
\[
I_{11}(x,t)\le3
\int_R^{t-x}d\beta
\int_\beta^{t+x}\frac{w(y,s)^{-p}}{(1+(\alpha-\beta)/2)^{1+a}}d\alpha.
\]
When $a>0$, the $\alpha$-integral is estimated as
\[
\int_\beta^{t+x}\frac{1}{(1+(\alpha-\beta)/2)^{1+a}}d\alpha
=\left[\frac{2}{-a}\left(1+\frac{\alpha-\beta}{2}\right)^{-a}\right]_{\alpha=\beta}^{\alpha=t+x}
\le\frac{2}{a},
\]
so that we have
\[
I_{11}(x,t)\le C(t-x-R)\le Cw(x,t)^{-1}E(T).
\]
When $a=0$, the $\alpha$-integral is estimated as
\[
\begin{array}{ll}
\d\int_\beta^{t+x}\frac{\log^p(\alpha+3R)}{1+(\alpha-\beta)/2}d\alpha
&\d\le\log^p(t+x+3R)\int_\beta^{t+x}\frac{1}{1+(\alpha-\beta)/2}d\alpha\\
&\d\le2\log^{p+1}(t+x+3R),
\end{array}
\]
so that we have
\[
I_{11}(x,t)\le C(t-x-R)\log^{p+1}(t+2R)\le Cw(x,t)^{-1}E(T).
\]
When $a<0$, the $\alpha$-integral is estimated as
\[
\begin{array}{l}
\d\int_\beta^{t+x}\frac{(\alpha+3R)^{-pa}}{(1+(\alpha-\beta)/2)^{1+a}}d\alpha\\
\d\le(t+x+3R)^{-pa}\left[\frac{2}{-a}\left(1+\frac{\alpha-\beta}{2}\right)^{-a}\right]_{\alpha=\beta}^{\alpha=t+x}\\
\d\le(t+x+3R)^{-pa}\cdot\frac{2^{1+a}}{-a}(t+x+2R)^{-a},
\end{array}
\]
so that we have
\[
I_{11}(x,t)\le C(t+2R)^{1-pa}(t+x+3R)^{-a}\le Cw(x,t)^{-1}E(T).
\]

\par
Next we shall deal with $I_{12}$ in $D_{\rm Int}$.
When $a>0$, we have
\[
I_{12}(x,t)=\int_{-R}^Rd\beta
\int_R^{t+x}\frac{1}{(1+(\alpha-\beta)/2)^{1+a}}d\alpha,
\]
so that the estimate is the same as $I_{01}$ in the proof of Proposition \ref{prop:apriori_0}
which implies that
\[
I_{12}(x,t)\le C\le Cw(x,t)^{-1}E(T).
\]
When $a=0$, we have
\[
I_{12}(x,t)=\int_{-R}^Rd\beta
\int_R^{t+x}\frac{\log^p(\alpha+3R)}{1+(\alpha-\beta)/2}d\alpha,
\]
so that
\[
I_{12}(x,t)\le2R\log^p(t+x+3R)
\int_R^{t+x}\frac{1}{1+(\alpha-R)/2}d\alpha
\]
follows, which implies
\[
I_{12}(x,t)\le C\log^{p+1}(t+x+3R)\le Cw(x,t)^{-1}E(T).
\]
When $a<0$, we have
\[
I_{12}(x,t)=\int_{-R}^Rd\beta
\int_R^{t+x}\frac{(\alpha+3R)^{-pa}}{(1+(\alpha-\beta)/2)^{1+a}}d\alpha,
\]
so that
\[
I_{12}(x,t)\le(t+x+3R)^{-pa}\int_{-R}^Rd\beta
\int_R^{t+x}\frac{1}{(1+(\alpha-\beta)/2)^{1+a}}d\alpha
\]
follows. The $\alpha$-integral is the same as $I_{01}$ in the proof of Proposition \ref{prop:apriori_0},
so that
\[
I_{12}(x,t)\le C(t+x+3R)^{-pa-a}\le Cw(x,t)^{-1}E(T)
\]
follows.

\par
Similarly to the above, we shall estimate $I_{13}$ in $D_{\rm Int}$.
When $a>0$, we have
\[
I_{13}(x,t)=\int_R^{t-x}d\beta
\int_{-R}^R\frac{1}{(1-(\alpha-\beta)/2)^{1+a}}d\alpha
\le2R\int_R^{t-x}d\beta
\]
so that
\[
I_{13}(x,t)\le C(t-x-R)\le Cw(x,t)^{-1}E(T)
\]
follows.
When $a=0$, we have
\[
\begin{array}{ll}
I_{13}(x,t)
&\d=\int_R^{t-x}d\beta
\int_{-R}^R\frac{\log^p(\alpha+3R)}{1-(\alpha-\beta)/2}d\alpha\\
&\d\le2R\log^p4R\int_R^{t-x}\frac{1}{1+(\beta-R)/2}d\beta
\end{array}
\]
which implies
\[
I_{13}(x,t)\le C\log(t-x+R)\le Cw(x,t)^{-1}E(T).
\]
When $a<0$, we have
\[
\begin{array}{ll}
I_{13}(x,t)
&\d=\int_R^{t-x}d\beta
\int_{-R}^R\frac{(\alpha+3R)^{-pa}}{(1-(\alpha-\beta)/2)^{1+a}}d\alpha\\
&\d\le\frac{2(4R)^{-pa}}{-a}\int_R^{t-x}\left(1+\frac{\beta+R}{2}\right)^{-a}d\beta
\end{array}
\]
which implies
\[
I_{13}(x,t)\le C(t-x+3R)^{1-a}\le Cw(x,t)^{-1}E(T).
\]
\par
It is easy to estimate $I_{14}$ in $D_{\rm Int}$.
Extending the domain of the integral, we have
\[
I_{14}(x,t)\le\int_{-R}^Rd\beta
\int_{-R}^R\frac{w(|y|,s)^{-p}}{(1+|\alpha-\beta|/2)^{1+a}}d\alpha
\]
which implies 
\[
I_{14}(x,t)\le C\le w(x,t)^{-1}E(T).
\]
Summing up all the estimates, we obtain
\[
I(x,t)\le w(x,t)^{-1}E(T)\quad\mbox{for}\ (x,t)\in D_{\rm Int}.
\]

\par
Let us step into the estimates in $D_{\rm Ext}$.
When $a>0$, we have
\[
I_{21}(x,t)\le\int_{-R}^Rd\beta
\int_R^{t+x}\frac{1}{(1+(\alpha-\beta)/2)^{1+a}}d\alpha,
\]
so that the estimate is completely the same as $I_{01}$ in the proof of Proposition \ref{prop:apriori_0}.
Hence we obtain
\[
I_{21}(x,t)\le C\le Cw(x,t)^{-1}E(T).
\]
When $a=0$, we have
\[
\begin{array}{ll}
I_{21}(x,t)
&\d\le\int_{-R}^Rd\beta\int_R^{t+x}\frac{\log^p(\alpha+3R)}{1+(\alpha-\beta)/2}d\alpha\\
&\d\le\log^p(t+x+3R)\int_{-R}^Rd\beta\int_R^{t+x}\frac{1}{1+(\alpha-\beta)/2}d\alpha.
\end{array}
\]
Hence, similarly to the above, we obtain
\[
I_{21}(x,t)\le C\log^{p+1}(t+x+3R)\le Cw(x,t)^{-1}E(T).
\]
When $a<0$, we have
\[
\begin{array}{ll}
I_{21}(x,t)
&\d\le\int_{-R}^Rd\beta\int_R^{t+x}\frac{(\alpha+3R)^{-pa}}{(1+(\alpha-\beta)/2)^{1+a}}d\alpha\\
&\d\le(t+x+3R)^{-pa}\int_{-R}^Rd\beta\int_R^{t+x}\frac{1}{(1+(\alpha-\beta)/2)^{1+a}}d\alpha.
\end{array}
\]
Hence we obtain
\[
I_{21}(x,t)\le C(t+x+3R)^{-pa-a}\le Cw(x,t)^{-1}E(T).
\]
Moreover, it is easy to estimate $I_{22}$ in $D_{\rm Ext}$.
Extending the domain of the integral, we have
\[
I_{22}(x,t)\le\int_{-R}^Rd\beta
\int_{-R}^R\frac{w(|y|,s)^{-p}}{(1+|\alpha-\beta|/2)^{1+a}}d\alpha
\]
which implies
\[
I_{22}(x,t)\le C\le Cw(x,t)^{-1}E(T).
\]
Summing up all the estimates, we obtain
\[
I(x,t)\le Cw(x,t)^{-1}E(T)\quad\mbox{for}\ (x,t)\in D_{\rm Ext}.
\]

\par
Finally we shall estimate $I_3$ in $D_{\rm Ori}$, but this is almost trivial because of
\[
I_3(x,t)\le\int_{-R}^Rd\beta
\int_{-R}^R\frac{w(|y|,s)^{-p}}{(1+|\alpha-\beta|/2)^{1+a}}d\alpha.
\]
Hence we obtain
\[
I_3(x,t)\le C\le Cw(x,t)^{-1}E(T)
\]
which implies
\[
I(x,t)\le C\le Cw(x,t)^{-1}E(T)\quad\mbox{for}\ (x,t)\in D_{\rm Ori}.
\]
Therefore (\ref{basic-estimate}) is established as desired.
This completes the proof of Proposition \ref{prop:apriori}.
\hfill$\Box$.

%%%%%%%%%%%%%%%%%%%%%%%%%%%%%%%%%%%%%%%%%%%%%%%%%%
%%%%%%%%%%%%%%%%%%%%% SECTION4 %%%%%%%%%%%%%%%%%%%%%%%
%%%%%%%%%%%%%%%%%%%%%%%%%%%%%%%%%%%%%%%%%%%%%%%%%%

\section{Proofs of Theorems \ref{thm:upper-bound_non-zero} and \ref{thm:upper-bound_zero}}
\par
In this section, we shall investigate the upper bounds of the lifespan.
As stated at the end of Section 2, the upper bounds of the lifespan are also determined
by point-wise estimates of the solution in the interior domain, $D_{\rm Int}$ in (\ref{Dint}).
In fact, it follows from  (\ref{supp_initial}) and (\ref{linear}) that
\[
u(x,t)=\frac{\e}{2}\int_{\R}g(x)dx+L_a(|u|^p)(x,t)
\quad\mbox{for}\ (x,t)\in D_{\rm Int}.
\]
Throughout this section, we assume that
\begin{equation}
\label{D}
(x,t)\in D:=D_{\rm Int}\cap\{x>0\}\cap\{t-x>R\}.
\end{equation}
Making use of (\ref{equivalent}) and
introducing the characteristic coordinate by (\ref{coordinate}), we have that
\begin{equation}
\label{first}
u(x,t)\ge C_0
\int_R^{t-x}d\beta\int_\beta^{t+x}\frac{|u(y,s)|^p}{(1+(\alpha-\beta)/2)^{1+a}}d\alpha
+J(x,t),
\end{equation}
where
\begin{equation}
\label{C_0}
C_0:=\frac{1}{8}\left(\frac{1}{\sqrt{2}}\right)^{\max\{0,-(1+a)\}}>0
\end{equation}
and
\begin{equation}
\label{J}
J(x,t):=
C_0\int_0^Rd\beta\int_\beta^{t+x}\frac{|u(y,s)|^p}{(1+(\alpha-\beta)/2)^{1+a}}d\alpha
+\frac{\e}{2}\int_{\R}g(x)dx.
\end{equation}
Employing this integral inequality, we shall estimate the lifespan from above.

%%% subsection 4.1 %%%

\subsection{Proof of Theorem \ref{thm:upper-bound_non-zero}}
Let $u=u(x,t)\in C^2(\R\times[0,T))$ be a solution of (\ref{IVP}).
It follows from (\ref{positive_non-zero}), (\ref{first}) and (\ref{J}) that
\begin{equation}
\label{frame1}
u(x,t)\ge C_0
\int_R^{t-x}d\beta\int_\beta^{t+x}\frac{|u(y,s)|^p}{(1+(\alpha-\beta)/2)^{1+a}}d\alpha+C_g\e
\end{equation}
for $(x,t)\in D$, where
\[
C_g:=\frac{1}{2}\int_{\R}g(x)dx>0.
\]

\vskip10pt
\par\noindent
{\bf Case 1.} $\v{a>0.}$
\par
Let
\begin{equation}
\label{D_R}
(x,t)\in D_R:=D\cap\{x\le R\}.
\end{equation}
Assume that an estimate
\begin{equation}
\label{step_n1}
u(x,t)\ge M_n\{(t-x-R)x\}^{a_n}
\quad\mbox{for}\ (x,t)\in D_R
\end{equation}
holds, where $a_n\ge0$ and $M_n>0$.
The sequences $\{a_n\}$ and $\{M_n\}$ are defined later. 
Then it follows from (\ref{frame1}) and (\ref{step_n1}) that
\[
u(x,t)\ge C_0M_n^p
\int_R^{t-x}(\beta-R)^{pa_n}d\beta
\int_\beta^{\beta+2x}\frac{\{(\alpha-\beta)/2\}^{pa_n}}{(1+(\alpha-\beta)/2)^{1+a}}d\alpha.
\]
Note that the domain of the integral is included in $D_R$, that is,
\[
\{(y,s):R\le \beta=s-y\le t-x,\ \beta\le \alpha=s+y\le \beta+2x\}\subset D_R
\]
for $(x,t)\in D_R$.
Since
\[
\int_\beta^{\beta+2x}\frac{\{(\alpha-\beta)/2\}^{pa_n}}{(1+(\alpha-\beta)/2)^{1+a}}d\alpha
\ge\frac{1}{(1+x)^{1+a}}\cdot\frac{2}{pa_n+1}x^{pa_n+1},
\]
we have
\[
u(x,t)\ge\frac{C_1M_n^p}{(pa_n+1)^2}\{(t-x-R)x\}^{pa_n+1}
\quad\mbox{for}\ (x,t)\in D_R,
\]
where
\begin{equation}
\label{C_1}
C_1:=\frac{2C_0}{(1+R)^{1+a}}>0.
\end{equation}
Therefore, if $\{a_n\}$ is defined by
\begin{equation}
\label{a_n}
a_{n+1}=pa_n+1,\ a_1=0,
\end{equation}
then (\ref{step_n1}) holds for all $n\in\N$ as far as $M_n$ satisfies
\begin{equation}
\label{M_n}
M_{n+1}\le\frac{C_1M_n^p}{(pa_n+1)^2}.
\end{equation}
In view of (\ref{frame1}), we note that (\ref{step_n1}) holds for $n=1$ with
\begin{equation}
\label{M_1}
M_1:=C_g\e.
\end{equation}

\par
Let us fix $\{M_n\}$.
It follows from (\ref{a_n}) that
\[
a_n=\frac{p^{n-1}-1}{p-1}\quad(n\in\N)
\]
which implies
\[
pa_n+1=a_{n+1}\le\frac{p^n}{p-1}.
\]
In view of (\ref{M_n}) and (\ref{M_1}), one of the choice of the definition of $\{M_n\}$ is
\begin{equation}
\label{def:M_n}
M_{n+1}=C_2p^{-2n}M_n^p,\ M_1=C_g\e,
\end{equation}
where
\begin{equation}
\label{C_2}
C_2:=(p-1)^2C_1>0.
\end{equation}
Hence we obtain that $M_n>0$ for all $n\in\N$ and
\[
\log M_{n+1}=\log C_2-2n\log p+p\log M_n
\]
which implies
\[
\begin{array}{ll}
\log M_{n+1}
&=(1+p+\cdots+p^{n-1})\log C_2\\
&\quad-2\{n+p(n-1)+\cdots+p^{n-1}(n-n+1)\}\log p+p^n\log M_1\\
&\d=\frac{p^n-1}{p-1}\log C_2-2p^n\log p\sum_{j=1}^n\frac{j}{p^j}+p^n\log M_1\\
&\d\ge-\frac{1}{p-1}\log C_2+p^n\left\{\frac{1}{p-1}\log C_2-2S_p\log p+\log M_1\right\},
\end{array}
\]
where
\begin{equation}
\label{S_p}
S_p:=\sum_{j=1}^\infty\frac{j}{p^j}<\infty
\end{equation}
because of d'Alembert's criterion.

\par
Therefore it follows from (\ref{step_n1}) that
\[
u(x,t)\ge C_3\{(t-x-R)x\}^{-1/(p-1)}\exp\left\{K_1(x,t)p^{n-1}\right\}
\quad\mbox{for}\ (x,t)\in D_R,
\]
where
\begin{equation}
\label{C_3}
C_3:=\exp\left(-\frac{1}{p-1}\log C_2\right)>0
\end{equation}
and
\begin{equation}
\label{K1}
\begin{array}{ll}
K_1(x,t):=&\d\frac{1}{p-1}\log\{(t-x-R)x\}\\
&\d+\frac{1}{p-1}\log C_2-2S_p\log p+\log(C_g\e).
\end{array}
\end{equation}
If there exists a point $(x_0,t_0)\in D_R$ such that
\[
K_1(x_0,t_0)>0,
\]
we have a contradiction
\[
u(x_0,t_0)=\infty
\]
by letting $n\rightarrow\infty$, so that $T<t_0$.
Let us set
\begin{equation}
\label{restriction}
x_0=R\quad\mbox{and}\quad t_0\ge4R.
\end{equation}
 Then $K_1(R,t_0)>0$ is equivalent to
\[
(t_0-2R)RC_2\exp\{-2(p-1)S_p\log p\}(C_g\e)^{p-1}>1.
\]
This condition follows from
\begin{equation}
\label{blow-up}
t_0>2R^{-1}C_2^{-1}\exp\{2(p-1)S_p\log p\}(C_g)^{1-p}\e^{-(p-1)}.
\end{equation}
We note that (\ref{blow-up}) is stronger than $t_0\ge4R$ for
\[
0<\e\le\e_3
\]
where $\e_3$ is defined by
\[
4R=2R^{-1}C_2^{-1}\exp\{2(p-1)S_p\log p\}(C_g)^{1-p}\e_3^{-(p-1)}.
\]
It is easy to see that $(R,t_0)\in D_R$ with $t_0$ satisfying (\ref{blow-up}).
The proof for $a>0$ is now completed.

\vskip10pt
\par\noindent
{\bf Case 2.} $\v{a=0.}$
\par
Assume that an estimate
\begin{equation}
\label{step_n2}
u(x,t)\ge M_n\{(t-x-R)\log(1+x)\}^{a_n}
\quad\mbox{for}\ (x,t)\in D
\end{equation}
holds, where $a_n\ge0$ and $M_n>0$.
The sequences $\{a_n\}$ and $\{M_n\}$ are defined later. 
Then it follows from (\ref{frame1}) and (\ref{step_n2}) that
\[
u(x,t)\ge C_0M_n^p
\int_R^{t-x}(\beta-R)^{pa_n}d\beta
\int_\beta^{t+x}\frac{\log^{pa_n}(1+(\alpha-\beta)/2)}{1+(\alpha-\beta)/2}d\alpha.
\]
Note that the domain of the integral is included in $D$, that is,
\[
\{(y,s):R\le \beta=s-y\le t-x,\ \beta\le \alpha=s+y\le t+x\}\subset D
\]
for $(x,t)\in D$.
Since
\[
\int_\beta^{t+x}\frac{\log^{pa_n}(1+(\alpha-\beta)/2)}{1+(\alpha-\beta)/2}d\alpha
=\frac{2}{pa_n+1}\log^{pa_n+1}\left(1+\frac{t+x-\beta}{2}\right),
\]
we have
\[
u(x,t)\ge\frac{2C_0M_n^p}{(pa_n+1)^2}\{(t-x-R)\log(1+x)\}^{pa_n+1}
\quad\mbox{for}\ (x,t)\in D.
\]

\par
Hence we can employ the same definitions of $\{M_n\}$ and $\{a_n\}$
as Case 1 in which $C_1$ is replaced with $2C_0$, so that we have
\[
u(x,t)\ge C_4\{(t-x-R)\log(1+x)\}^{-1/(p-1)}\exp\left\{K_2(x,t)p^{n-1}\right\}
\quad\mbox{for}\ (x,t)\in D,
\]
where
\begin{equation}
\label{C_4C_5}
C_4:=\exp\left(-\frac{1}{p-1}\log C_5\right)>0,\ C_5:=2(p-1)^2C_0>0
\end{equation}
and
\[
\begin{array}{ll}
K_2(x,t):=
&\d\frac{1}{p-1}\log\{(t-x-R)\log(1+x)\}\\
&\d+\frac{1}{p-1}\log C_5-2S_p\log p+\log(C_g\e),
\end{array}
\]
where $S_p$ is the one in (\ref{S_p}).

\par
Therefore the same argument as Case 1 is valid.
The difference appears only in finding $(x_0,t_0)\in D$ with $K_2(x_0,t_0)>0$.
Let
\[
t_0=2x_0\quad\mbox{and}\quad t_0\ge4R.
\]
Then, since we have
\[
(t_0-x_0-R)\log(1+x_0)\ge\frac{t_0}{4}\log\left(1+\frac{t_0}{2}\right)\ge\frac{t_0}{8}\log(2+t_0),
\]
$K_2(t_0/2,t_0)>0$ follows from
\[
\phi(t_0)=t_0\log(2+t_0)>8C_5^{-1}\exp\{2(p-1)S_p\log p\}(C_g)^{1-p}\e^{-(p-1)}.
\]
This completes the proof for $a=0$.

\vskip10pt
\par\noindent
{\bf Case 3.} $\v{a<0.}$
\par
This case is almost similar to Case 2.
Assume that an estimate
\begin{equation}
\label{step_n3}
u(x,t)\ge M_n\left\{(t-x-R)\frac{x^{1-a}}{1+t+x}\right\}^{a_n}
\quad\mbox{for}\ (x,t)\in D
\end{equation}
holds, where $a_n\ge0$ and $M_n>0$.
The sequences $\{a_n\}$ and $\{M_n\}$ are defined later. 
Then it follows from (\ref{frame1}) and (\ref{step_n3}) that
\[
u(x,t)\ge C_0M_n^p
\int_R^{t-x}(\beta-R)^{pa_n}d\beta
\int_\beta^{t+x}\frac{\{(\alpha-\beta)/2\}^{(1-a)pa_n}}{(1+(\alpha-\beta)/2)^{1+a}(1+\alpha)^{pa_n}}d\alpha.
\]
Note that the domain of the integral is included in $D$, that is,
\[
\{(y,s):R\le \beta=s-y\le t-x,\ \beta\le \alpha=s+y\le t+x\}\subset D
\]
for $(x,t)\in D$.
Since
\[
\begin{array}{l}
\d\int_\beta^{t+x}\frac{\{(\alpha-\beta)/2\}^{(1-a)pa_n}}{(1+(\alpha-\beta)/2)^{1+a}(1+\alpha)^{pa_n}}d\alpha\\
\d\ge\frac{1}{(1+t+x)^{pa_n+1}}
\int_\beta^{t+x}\left(\frac{\alpha-\beta}{2}\right)^{-a+(1-a)pa_n}d\alpha\\
\d=\frac{2}{(1-a)(pa_n+1)(1+t+x)^{pa_n+1}}\left(\frac{t+x-\beta}{2}\right)^{(1-a)(pa_n+1)}
\end{array}
\]
hold, we have
\[
u(x,t)\ge\frac{2C_0M_n^p}{(1-a)(pa_n+1)^2}\left\{(t-x-R)\frac{x^{1-a}}{1+t+x}\right\}^{pa_n+1}
\quad\mbox{for}\ (x,t)\in D.
\]

\par
Hence we can employ the same definitions of $\{M_n\}$ and $\{a_n\}$
as Case 1 in which $C_1$ is replaced with $2C_0/(1-a)$, so that we have
\[
u(x,t)\ge C_6\left\{(t-x-R)\frac{x^{1-a}}{1+t+x}\right\}^{-1/(p-1)}\exp\left\{K_3(x,t)p^{n-1}\right\}
\quad\mbox{for}\ (x,t)\in D,
\]
where
\begin{equation}
\label{C_6C_7}
C_6:=\exp\left(-\frac{1}{p-1}\log C_7\right)>0,\ C_7:=\frac{2(p-1)^2C_0}{1-a}>0
\end{equation}
and
\[
\begin{array}{ll}
K_3(x,t):=
&\d\frac{1}{p-1}\log\left\{(t-x-R)\frac{x^{1-a}}{1+t+x}\right\}\\
&\d+\frac{1}{p-1}\log C_7-2S_p\log p+\log(C_g\e),
\end{array}
\]
where $S_p$ is the one in (\ref{S_p}).

\par
Therefore the same argument as Case 1 is valid.
The difference appears only in finding $(x_0,t_0)\in D$ with $K_3(x_0,t_0)>0$.
Let
\[
t_0=2x_0\quad\mbox{and}\quad t_0\ge4R.
\]
Then, since we have
\[
(t_0-x_0-R)\frac{x_0^{1-a}}{1+t_0+x_0}
\ge\frac{t_0}{4}\left(\frac{t_0}{2}\right)^{1-a}\frac{1}{R+t_0+t_0+R},
\]
$K_3(t_0/2,t_0)>0$ follows from
\[
t_0^{1-a}>5\cdot2^{2-a}C_7^{-1}\exp\{2(p-1)S_p\log p\}(C_g)^{1-p}\e^{-(p-1)}.
\]
This completes the proof for $a<0$.
\hfill$\Box$

%%% subsection 4.2 %%%
\subsection{Proof of Theorem \ref{thm:upper-bound_zero}}
The proof is almost similar to the one of Theorem \ref{thm:upper-bound_non-zero}.
Let $u=u(x,t)\in C^2(\R\times[0,T))$ be a solution of (\ref{IVP}).
Since the assumption on the initial data in (\ref{positive_zero}) yields
\begin{equation}
\label{lower-bound_linear}
u^0(x,t)=\frac{1}{2}\{f(x+t)+f(x-t)\}\ge\frac{1}{2}f(x-t)
\quad\mbox{for}\quad(x,t)\in\R\times[0,\infty),
\end{equation}
it follows from (\ref{first}) and (\ref{J}) that
\[
u(x,t)\ge\e u^0(x,t)\ge\frac{\e}{2}f(x-t)\quad\mbox{for}\quad(x,t)\in\R\times[0,T)
\]
and
\begin{equation}
\label{frame2}
u(x,t)\ge C_0
\int_R^{t-x}d\beta\int_\beta^{t+x}\frac{|u(y,s)|^p}{(1+(\alpha-\beta)/2)^{1+a}}d\alpha
+\frac{C_0}{2^p}\e^pJ'(x,t)
\end{equation}
for $(x,t)\in D$, where $D$, $C_0$ are defined in (\ref{D}), (\ref{C_0}) respectively, and
\begin{equation}
\label{J'}
J'(x,t):=\int_0^Rf(-\beta)^pd\beta\int_\beta^{t+x}\frac{1}{(1+(\alpha-\beta)/2)^{1+a}}d\alpha.
\end{equation}
Note that, without loss of the generality, we may assume that
\begin{equation}
\label{positive_zero_add}
f(x)\not\equiv0\quad\mbox{for}\ x\in(-R,0).
\end{equation}
Because, if not, we have to assume that
\[
f(x)\not\equiv0\quad\mbox{for}\ x\in(0,R).
\]
Therefore we obtain all the estimates below for $x<0$ by replacing $x$ with $-x$.
Because, taking $f(x+t)$ instead of $f(x-t)$ in (\ref{lower-bound_linear}),
we have, in stead of (\ref{frame2}), that
\[
u(x,t)\ge C_0
\int_R^{t+x}d\alpha\int_\alpha^{t-x}\frac{|u(y,s)|^p}{(1+(\beta-\alpha)/2)^{1+a}}d\beta
+\frac{C_0}{2^p}\e^pJ''(x,t),
\]
where
\[
J''(x,t):=\int_0^Rf(\alpha)^pd\alpha\int_\alpha^{t-x}\frac{1}{(1+(\beta-\alpha)/2)^{1+a}}d\alpha.
\]
This implies the symmetry of the domain as well as the estimates.

\vskip10pt
\par\noindent
{\bf Case 1.} $\v{a>0.}$
\par
In this case, we assume (\ref{D_R}) again.
Recall that
\[
t+x=t-x+2x\ge R+2x\ge\beta+2x
\quad\mbox{for}\ (x,t)\in D_R\ \mbox{and}\ \beta\in[0,R].
\]
Since
\[
\begin{array}{ll}
\d\int_\beta^{t+x}\frac{1}{(1+(\alpha-\beta)/2)^{1+a}}d\alpha
&\d\ge\int_\beta^{\beta+2x}\frac{1}{(1+(\alpha-\beta)/2)^{1+a}}d\alpha\\
&\d\ge\frac{2x}{(1+x)^{1+a}},
\end{array}
\]
holds for $\beta\in[0,R]$, it follows from (\ref{frame2}) and (\ref{positive_zero_add}) that
\begin{equation}
\label{frame21}
u(x,t)\ge C_0
\int_R^{t-x}d\beta\int_\beta^{\beta+2x}\frac{|u(y,s)|^p}{(1+(\alpha-\beta)/2)^{1+a}}d\alpha
+C_f\e^px
\end{equation}
for $(x,t)\in D_R$, where
\[
C_f:=\frac{2C_0}{2^p(1+R)^{1+a}}\int_0^Rf(-\beta)^pd\beta>0
\]
because of
\[
t+x=t-x+2x\ge\beta+2x
\quad\mbox{for}\ (x,t)\in D_R\ \mbox{and}\ \beta\in[R,t-x].
\]

\par
From now on, we employ the same argument as Case 1
of the proof of Theorem \ref{thm:upper-bound_zero}.
Instead of (\ref{step_n1}), assume that an estimate
\begin{equation}
\label{step_n4}
u(x,t)\ge M_n(t-x-R)^{a_n}x^{b_n}
\quad\mbox{for}\ (x,t)\in D_R
\end{equation}
holds, where $a_n\ge0$, $b_n>0$ and $M_n>0$.
The sequences $\{a_n\},\{b_n\}$ and $\{M_n\}$ are defined later. 
Then it follows from (\ref{frame21}), (\ref{step_n4}) and the same computations
as Case 1 of the proof of Theorem \ref{thm:upper-bound_zero} that
(\ref{step_n4}) holds for all $n\in\N$ provided
\[
\left\{
\begin{array}{l}
a_{n+1}=pa_n+1,\ a_1=0,\\
b_{n+1}=pb_n+1,\ b_1=1
\end{array}
\right.
\]
and
\[
M_{n+1}\le\frac{C_1M_n^p}{(pa_n+1)(pb_n+1)},\ M_1=C_f\e^p,
\]
where $C_1$ is the one in (\ref{C_1}).
It is easy to see that
\[
a_n=\frac{p^{n-1}-1}{p-1},\ b_n=\frac{p^n-1}{p-1}\quad (n\in\N)
\]
which implies
\[
(pa_n+1)(pb_n+1)\le(pb_n+1)^2=b_{n+1}^2\le\frac{p^{2(n+1)}}{(p-1)^2}.
\]
Hence $M_n$ in this case should be defined by
\[
M_{n+1}=C_2p^{-2(n+1)}M_n^p,\ M_1=C_f\e^p,
\]
where $C_2$ is the one in (\ref{C_2}).

\par
Therefore it follows from (\ref{step_n4}) that
\[
u(x,t)\ge C_3\{(t-x-R)x\}^{-1/(p-1)}\exp\left\{K_4(x,t)p^{n-1}\right\}
\quad\mbox{for}\ (x,t)\in D_R,
\]
where $C_3$ is the one in (\ref{C_3}) and
\[
\begin{array}{ll}
K_4(x,t):=&\d\frac{1}{p-1}\log\{(t-x-R)x^p\}\\
&\d+\frac{1}{p-1}\log C_2-2S_p'\log p+\log(C_f\e^p),
\end{array}
\]
where
\begin{equation}
\label{S_p'}
S_p':=\sum_{j=1}^\infty\frac{j+1}{p^j}<\infty.
\end{equation}
Assume (\ref{restriction}) again.
Then, $K_4(R,t_0)>0$ is equivalent to
\[
(t_0-2R)R^pC_2\exp\{-2(p-1)S_p'\log p\}(C_f\e^p)^{p-1}>1.
\]
This condition follows from
\begin{equation}
\label{blow-up2}
t_0>2R^{-p}C_2^{-1}\exp\{2(p-1)S_p'\log p\}(C_f)^{1-p}\e^{-p(p-1)}.
\end{equation}
We note that this is stronger than $t_0\ge4R$ for
\[
0<\e\le\e_4
\]
where $\e_4$ is defined by
\[
4R=2R^{-p}C_2^{-1}\exp\{2(p-1)S_p'\log p\}(C_f)^{1-p}\e_4^{-p(p-1)}.
\]
It is easy to see that $(R,t_0)\in D_R$ with $t_0$ satisfying (\ref{blow-up2}).
The proof for $a>0$ is now completed.

\vskip10pt
\par\noindent
{\bf Case 2.} $\v{a=0.}$
\par
Recall that
\begin{equation}
\label{relation_D}
t+x-\beta=t-x-R+2x\ge2x
\quad\mbox{for}\ (x,t)\in D\ \mbox{and}\ \beta\in[0,R].
\end{equation}
Since
\[
\int_\beta^{t+x}\frac{1}{1+(\alpha-\beta)/2}d\alpha
\ge2\log\left(1+\frac{t+x-\beta}{2}\right)\ge2\log(1+x)
\]
holds for $(x,t)\in D$ and $\beta\in[0,R]$,
it follows from (\ref{frame2}) and (\ref{positive_zero_add}) that
\begin{equation}
\label{frame22}
u(x,t)\ge C_0
\int_R^{t-x}d\beta\int_\beta^{t+x}\frac{|u(y,s)|^p}{(1+(\alpha-\beta)/2)^{1+a}}d\alpha
+C_f'\e^p\log(1+x)
\end{equation}
for $(x,t)\in D$, where
\[
C_f':=\frac{2C_0}{2^p}\int_0^Rf(-\beta)^pd\beta>0.
\]

\par
From now on, we employ the same argument as Case 2
of the proof of Theorem \ref{thm:upper-bound_zero}.
Instead of (\ref{step_n2}), assume that an estimate
\begin{equation}
\label{step_n5}
u(x,t)\ge M_n(t-x-R)^{a_n}\log^{b_n}(1+x)
\quad\mbox{for}\ (x,t)\in D
\end{equation}
holds, where $a_n\ge0,\ b_n>0$ and $M_n>0$.
The sequences $\{a_n\},\{b_n\}$ and $\{M_n\}$ are defined later. 
Then it follows from (\ref{frame22}), (\ref{step_n5}) and the same computations
as Case 2 of the proof of Theorem \ref{thm:upper-bound_zero} that
(\ref{step_n5}) holds for all $n\in\N$ provided
\[
\left\{
\begin{array}{l}
a_{n+1}=pa_n+1,\ a_1=0,\\
b_{n+1}=pb_n+1,\ b_1=1
\end{array}
\right.
\]
and
\[
M_{n+1}\le\frac{2C_0M_n^p}{(pa_n+1)(pb_n+1)},\ M_1=C_f'\e^p.
\]

\par
Hence we can employ the same definition of $\{M_n\}$
as Case 1 in which $C_1$ is replaced with $2C_0$, so that we have
\[
u(x,t)\ge C_4\{(t-x-R)\log(1+x)\}^{-1/(p-1)}\exp\left\{K_5(x,t)p^{n-1}\right\}
\quad\mbox{for}\ (x,t)\in D,
\]
where $C_4,C_5$ are defined in (\ref{C_4C_5}) and
\[
\begin{array}{ll}
K_5(x,t):=
&\d\frac{1}{p-1}\log\{(t-x-R)\log^p(1+x)\}\\
&\d+\frac{1}{p-1}\log C_5-2S_p'\log p+\log(C_f'\e^p),
\end{array}
\]
where $S_p'$ is the one in (\ref{S_p'}).

\par
Therefore the same argument as Case 1 is valid.
The difference appears only in finding $(x_0,t_0)\in D$ with $K_5(x_0,t_0)>0$.
Let
\[
t_0=2x_0\quad\mbox{and}\quad t_0\ge4R.
\]
Then, since we have
\[
(t_0-x_0-R)\log^p(1+x_0)\ge\frac{t_0}{4}\log^p\left(1+\frac{t_0}{2}\right)\ge\frac{t_0}{4\cdot2^p}\log(2+t_0),
\]
$K_2(t_0/2,t_0)>0$ follows from
\[
\psi_p(t_0)=t_0\log^p(2+t_0)>4\cdot2^pC_5^{-1}\exp\{2(p-1)S_p'\log p\}(C_f')^{1-p}\e^{-p(p-1)}.
\]
This completes the proof for $a=0$.

\vskip10pt
\par\noindent
{\bf Case 3.} $\v{a<0.}$
\par
This case is almost similar to Case 2.
Since
\[
\begin{array}{ll}
\d\int_\beta^{t+x}\frac{1}{(1+(\alpha-\beta)/2)^{1+a}}d\alpha
&\d\ge\frac{1}{1+t+x}\int_\beta^{t+x}\left(\frac{\alpha-\beta}{2}\right)^{-a}d\alpha\\
&\d=\frac{1}{1+t+x}\cdot\frac{2}{1-a}\left(\frac{t+x-\beta}{2}\right)^{1-a}\\
&\d\ge\frac{2}{1-a}\cdot\frac{x^{1-a}}{1+t+x}
\end{array}
\]
holds for $(x,t)\in D$ and $\beta\in[0,R]$ because of (\ref{relation_D}),
it follows from (\ref{frame2}) and (\ref{positive_zero_add}) that
\begin{equation}
\label{frame23}
u(x,t)\ge C_0
\int_R^{t-x}d\beta\int_\beta^{t+x}\frac{|u(y,s)|^p}{(1+(\alpha-\beta)/2)^{1+a}}d\alpha
+C_f''\e^p\frac{x^{1-a}}{1+t+x}
\end{equation}
for $(x,t)\in D$, where
\[
C_f'':=\frac{2C_0}{2^p(1-a)}\int_0^Rf(-\beta)^pd\beta>0.
\]

\par
From now on, we employ the same argument as Case 3
of the proof of Theorem \ref{thm:upper-bound_zero}.
Instead of (\ref{step_n3}), assume that an estimate
\begin{equation}
\label{step_n6}
u(x,t)\ge M_n(t-x-R)^{a_n}\left(\frac{x^{1-a}}{1+t+x}\right)^{b_n}
\quad\mbox{for}\ (x,t)\in D
\end{equation}
holds, where $a_n\ge0,\ b_n>0$ and $M_n>0$.
The sequences $\{a_n\},\{b_n\}$ and $\{M_n\}$ are defined later. 
Then it follows from (\ref{frame23}), (\ref{step_n6}) and the same computations
as Case 3 of the proof of Theorem \ref{thm:upper-bound_zero} that
(\ref{step_n6}) holds for all $n\in\N$ provided
\[
\left\{
\begin{array}{l}
a_{n+1}=pa_n+1,\ a_1=0,\\
b_{n+1}=pb_n+1,\ b_1=1
\end{array}
\right.
\]
and
\[
M_{n+1}\le\frac{2C_0M_n^p}{(1-a)(pa_n+1)(pb_n+1)},\ M_1=C_f''\e^p.
\]

\par
Hence we can employ the same definition of $\{M_n\}$
as Case 3 in which $C_1$ is replaced with $2C_0$, so that we have
\[
u(x,t)\ge C_6\left\{(t-x-R)\frac{x^{1-a}}{1+t+x}\right\}^{-1/(p-1)}\exp\left\{K_6(x,t)p^{n-1}\right\}
\quad\mbox{for}\ (x,t)\in D,
\]
where $C_6,C_7$ are defined in (\ref{C_6C_7}) and
\[
\begin{array}{ll}
K_6(x,t):=
&\d\frac{1}{p-1}\log\left\{(t-x-R)\left(\frac{x^{1-a}}{1+t+x}\right)^p\right\}\\
&\d+\frac{1}{p-1}\log C_7-2S_p'\log p+\log(C_f''\e^p),
\end{array}
\]where $S_p'$ is the one in (\ref{S_p'}).

\par
Therefore the same argument as Case 1 and 2 is valid.
The difference appears only in finding $(x_0,t_0)\in D$ with $K_6(x_0,t_0)>0$.
Let
\[
t_0=2x_0\quad\mbox{and}\quad t_0\ge4R.
\]
Then, since we have
\[
(t_0-x_0-R)\left(\frac{x_0^{1-a}}{1+t_0+x_0}\right)^p
\ge\frac{t_0}{4}\left(\frac{t_0}{2}\right)^{p(1-a)}\frac{1}{(R+t_0+t_0+R)^p},
\]
$K_6(t_0/2,t_0)>0$ follows from
\[
t_0^{1-a}>5^p\cdot2^{2-pa}C_7^{-1}\exp\{2(p-1)S_p'\log p\}(C_f'')^{1-p}\e^{-p(p-1)}.
\]
This completes the proof for $a<0$.
\hfill$\Box$

%%%%%%%%%%%%%i%%%%%%%%%%%%%%%%%%%%%%%
%%%%%%%%%%%% Acknowledgement %%%%%%%%%%%%%%%
%%%%%%%%%%%%%%%%%%%%%%%%%%%%%%%%%%%%%
\section*{Acknowledgement}
\par
The third author is partially supported
by the Grant-in-Aid for Scientific Research (B) (No.18H01132), 
Japan Society for the Promotion of Science.
All the authors thank to the referee for pointing out many typos.

%%%%%%%%%%%%%%%%%%%%%%%%%%%%%%%%%%%%%%
%%%%%%%%%%%% References %%%%%%%%%%%%%%%%%%%%
%%%%%%%%%%%%%%%%%%%%%%%%%%%%%%%%%%%%%%

\bibliographystyle{plain}

\end{document}